\documentclass{amsart}
\usepackage{amssymb,amsmath}
\input{diagrams.tex}
\diagramstyle[scriptlabels]
\def \bbbc{{\mathbb C}}
\def \bbbp{{\mathbb P}}
\def \bbbs{{\mathbb S}}
\def\cI{{\mathcal I}}
\def\cL{{\mathcal L}}
\def\cM{{\mathcal M}}
\def\cO{{\mathcal O}}
\def\fm{{\mathfrak m}}

\def\sbullet{{\scriptscriptstyle \bullet}}
\def\Sing{\mathop{\rm Sing}\nolimits}
\def\Gr{\mathop{\rm Gr}\nolimits}
\def\codim{\mathop{\rm codim}\nolimits}
\def\red{{\mathop{\rm red}\nolimits}}
\def\Reg{\mathop{\rm Reg}\nolimits}
\def\reg{{\mathop{\rm reg}\nolimits}}
\def\Hom{\mathop{\rm Hom}\nolimits}
\def\Res{\mathop{\rm Res}\nolimits}
\def\Ext{\mathop{\rm Ext}\nolimits}
\def\image{\mathop{\rm Im}\nolimits}
\def\depth{\mathop{\rm depth}\nolimits}
\def\Projan{\mathop{\rm Projan}\nolimits}
\def\hto{\hookrightarrow}
\def\lto{\longrightarrow}

\theoremstyle{plain}

\newtheorem{theorem}{Theorem}[section]

\newtheorem{lemma}[theorem]{Lemma}
\newtheorem{proposition}[theorem]{Proposition}

\theoremstyle{definition}

\newtheorem{remark}[theorem]{Remark}

\newtheorem{sit}[theorem]{}

\begin{document}

\title{Contact Singularities}

\author{Fr\'ed\'eric Campana\and Hubert Flenner}
\thanks{The authors were supported by the Schwerpunktprogramm ``Globale
Methoden der komplexen Analysis" of the DFG.} 

\address{Fr\'ed\'eric Campana:
D\'epartement de
Math\'ematiques, Universit\'e Nancy 1, BP 239, 54506
Vandoeuvre-les-Nancy C\'edex, France}
\email{campana@iecn.u-nancy.fr}
\address{Hubert Flenner: Fakult\"at f\"ur Mathematik der Ruhr-Universit\"at,
Universit\"atsstr.\ 150, Geb.\ NA 2/72, 44780 Bochum, Germany}
\email{Hubert.Flenner@ruhr-uni-bochum.de}

\maketitle

\begin{abstract} 
A contact singularity is a normal singularity $(V,0)$
together with a holomorphic contact form $\eta$ on $V\backslash \Sing
V$ in a neighbourhood of $0$, i.e.\ $\eta\wedge (d\eta)^r$ has no
zero, where $\dim V=2r+1$. The main result of this paper is that
there are no isolated contact singularities. 
\end{abstract}

\section{Introduction}

Let $(V,0)$ be a normal complex singularity of dimension $n:=2r+1$. In
analogy with the notion of symplectic singularities, which were
studied by Beauville \cite{Bea0,Bea}, we call
$(V,0)$ a {\em contact singularity}\/ if there is a holomorphic
1-form on $V\backslash \Sing V$ in a neighbourhood of $0$ such that
$\eta\wedge (d\eta)^r$ has no zero on $V\backslash \Sing V$.
The main result of this paper is that {\em there are no isolated
contact singularities}, see \ref{c.4}. 

The key tool in proving this is a Hodge theoretic result for
rational singularities which may be of independent interest. Let
$\pi:X\to V$ be a resolution of an isolated rational
singularity $(V,0)$ such that the exceptional
set $E:=\pi^{-1}(0)$ is an SNC (=simple normal crossing) divisor. It
was shown by Steenbrink and van Straaten \cite{SSt} that then every
holomorphic $p$-form, say, $\eta$ on
$V\backslash \{0\}$ extends to a holomorphic form on the resolution
$X$. In \ref{1.1} we will show more strongly that for $p\ge 1$ even
$\pi^*(\eta)|E=0$ so that $\pi^*(\eta)$ belongs to $\Omega^p_X\langle
E\rangle(-E)$, where
$\Omega^p_X\langle E\rangle$ denotes the sheaf of meromorphic
$p$-forms with at most logarithmic poles along $E$ which are
holomorphic on
$X\backslash E$. This result should be viewed as a generalization of
the fact that on a Fano manifold there are no non-zero holomorphic
$p$-forms, $p\ge 1$, see Sect.\ 2.

We add a few remarks about the contents of the various sections.
Sect.\ 2 contains the aforementioned result about the extendability of
differential forms on rational singularities. In
Section 3 we study contact singularities. With the same arguments as
in the case of symplectic singularities \cite{Bea0,Bea} they are
Gorenstein and rational if $\codim\Sing V\ge 3$, see \ref{c.1}. We
treat the case of complete intersections more closely and show then
the main result described above. 

Sect.\ 4 contains a technical result on holomorphic forms for a
special class of isolated singularities which is used in Sect.\ 5 to
compute the number of moduli of an isolated symplectic singularity
that is resolved after one blowing up. It turns out that in this
case the number of moduli is a topological invariant given by
$\beta_2-1$, where $\beta_2$ is the second Betti number of the
exceptional set of the blowing up.

\section{Differential forms on rational singularities}

To motivate the main result of this section, let us first consider
the situation that $(V,0)\subseteq (\bbbc^N,0)$ is an isolated
singularity of dimension
$n\ge 2$ that is given by homogeneous equations. Blowing  up the point
$0\in V$ gives a resolution of singularities $\pi:X\to V$ with
exceptional set, say, $E$, which is a
projective manifold of dimension $n-1$. If $(V,0)$ is rational and
Gorenstein then by a result of Kempf \cite[Prop.\ on p.\
50]{KKMS} a nowhere vanishing holomorphic
$n$-form $\omega$ on
$V\backslash\{0\}$ extends to a holomorphic $n$-form on $X$ which
shows that
$\omega_X\cong\cO_X(kE)$ for some $k\ge 0$. By adjunction
$$
\omega_E\cong \omega_X(E)\cong\cO_X((k+1)E)
$$
and so the sheaf $\omega_E^{-1}$ is ample. Hence 
$E$ is a Fano manifold. In particular $E$ admits no global
holomorphic $p$-forms, $p\ge 1$. If $\eta$ is a
holomorphic $p$-form on $V\backslash\{0\}$ then by the results of
\cite{SSt} this form extends to a holomorphic $p$-form $\pi^*(\eta)$
on $X$. As $E$ is Fano, this form vanishes when restricted to $E$ and
so 
$$
\pi^*(\eta)\in\Gamma(X,\Omega^p_X\langle E\rangle(-E)).\leqno (*)
$$
Our aim is to generalize this observation to non-homogeneous
singularities. We will use the following notations.

Let $(V,0)$ be an arbitrary isolated singularity
and  consider as before a resolution of singularities $\pi:X\to V$
such that the exceptional set $E:=\pi^{-1}(0)$ is an SNC
divisor. 

The main result of this section is the following theorem.

\begin{theorem}\label{1.1}
Assume that $(V,0)$ is a rational isolated singularity of dimension
$n\ge 2$ and let $\eta$ be a holomorphic $p$-form on
$U:=V\backslash\{0\}$, where $p\ge 0$. Then
$\pi^*(\eta)$ extends to a holomorphic $p$-form in
$\Omega^p_X\langle E\rangle(-E)$.
\end{theorem}

Thus, if $\tilde \pi: Y\to V$ is any resolution of singularities and 
$F:=\tilde\pi^{-1}(0)$ then $\tilde\pi^*(\eta)$ extends to a
holomorphic form on $Y$ such the restriction $\tilde\pi^*(\eta)|F$
vanishes as a section of $\Omega^1_{F_\red}/Torsion$. This follows
easily by dominating $Y$ by a resolution $\pi:X\to V$ for which
$E=\pi^{-1}(0)$ is a simple normal crossing divisor.

Before proving the main theorem we remind the reader of some
basic results concerning mixed Hodge structures of isolated
singularities; for a general reference of the theory of mixed Hodge
structures see \cite{Elz}. 

Let $V$ be a
contractible Stein space and assume that $0\in V$ is the only
singular point of $V$. Let $\pi:X\to V$ be as before a resolution
of singularities such that $E:=\pi^{-1}(0)$ is an SNC divisor. By
\cite[1.9]{Ste1} (cf.\ also \cite[Sect.\ 5]{Elz}) the local
cohomology group $H^p_{\{0\}}(V,\bbbc)$  carries a mixed
Hodge structure. As $V$ is contractible this local cohomology group
is isomorphic to $H^{p-1}(U, \bbbc)$ for $p\ne 1$ whereas
$H^1_{\{0\}}(V,\bbbc)\cong H^{0}(U, \bbbc)/\bbbc$. Thus $H^{p}(U,
\bbbc)$ also carries a mixed Hodge structure. We have the following
facts.

\begin{theorem}\label{1.2}
{\em (1) (\cite[Sect.\ 1]{Ste1})} The spaces $H^p(E,\bbbc)$ and
$H^p(U,\bbbc)$ carry natural mixed Hodge structures, and the Hodge
filtrations are given by 
$$
F^pH^p(E,\bbbc)=H^0(E,\Omega^p_X/\Omega_X^p\langle E\rangle (-E))
\leqno (a)
$$
and 
$$
F^pH^p(U,\bbbc)=H^0(E,\Omega^p_X\langle E\rangle/\Omega_X^p\langle
E\rangle (-E))
\leqno (b)
$$
Moreover, the natural map 
$$
H^p(E,\bbbc)\cong H^p(X,\bbbc)\to
H^p(U,\bbbc) \leqno (c)
$$
is a morphism of mixed Hodge structures.

{\em (2) (\cite[1.3 and 1.4]{SSt})} If $\eta$ is a holomorphic
$p$-form on $U$ and if $p\le n-2$ then $\pi^*(\eta)$ extends to a
holomorphic $p$-form on $X$. If moreover $(V,0)$ has a rational
singularity then the same is true for $p=n-1$ and $p=n$.
\end{theorem}

\proof
 (1) was stated in this form in \cite[proof of
1.3]{SSt}. We add a few remarks how to deduce it from \cite{Ste1}.
(a) and (c) are immediate consequences of \cite[1.5 and 1.12]{Ste1},
whereas (b) is an easy consequence of the description of the
Hodge filtration given in \cite[1.9]{Ste1}. Note 
that the complex
$A^\sbullet_{\{0\}}(V)$ constructed in [loc.cit., 1.8] does not
compute the local cohomology $H^\sbullet_{\{0\}}(V,\bbbc)$ (as was
claimed there) but
$H^{\sbullet-1}(U,\bbbc)$ (the error comes from the isomorphism in
the second paragraph of 1.8 which is only valid for $k\ne 1$). 
\qed

Thus, if $\eta$ is as in \ref{1.1}, then the form $\pi^*(\eta)$
extends to a holomorphic form on $X$ by \ref{1.2} (2). Hence
\ref{1.1} is a consequence of the following lemma.

\begin{lemma}\label{1.3}
Under the assumptions of \ref{1.1} we have that
$$
H^0(X,\Omega^p_X)\cong H^0(X, \Omega^p_X\langle
E\rangle(-E)).
$$
\end{lemma}

\proof
Consider the $W$-filtration associated to the mixed Hodge structure
on the cohomology $H^*(E,\bbbc)$. Using the standard description of
the $W$-filtration (see e.g.\ \cite[3.5]{Elz}) it follows that
$$
\Gr_p^W H^p(E,\bbbc)=\image\left(
H^p(E,\bbbc)\lto \bigoplus_i H^p(E_i,\bbbc)\right)\,,
$$
where the sum is taken over all irreducible components $E_i$ of $E$.
This module carries a Hodge structure of weight $p$. In a first step
we will show that 
$$
F^p \Gr_p^W H^p(E,\bbbc)=0,\leqno (1)
$$
where $F^p \Gr_p^W H^p(E,\bbbc)$ is the $p$-th piece of the Hodge
filtration. By Hodge symmetry this piece is
isomorphic to
$$
M:=\frac{F^0 \Gr_p^W H^p(E,\bbbc)}{F^1 \Gr_p^W H^{p}(E,\bbbc)}.
$$
Using the fact that a morphism of Hodge structures is strict (cf.\
\cite[1.3.7 (iii)]{Elz}) the latter module is equal to the image of
$\Gr_p^W H^p(E,\bbbc)$ under the map $\beta$ in the
diagram
\begin{diagram}[s=8mm]
H^p(X,\bbbc)& \rTo^{\sim} &
H^p(E,\bbbc)&\rTo &\bigoplus H^p(E_i,\bbbc)\\
\dTo>{can} && \dTo>{\alpha=can} && \dTo>{\beta=can}\\
H^p(X,\cO_X)& \rTo&
H^p(E,\cO_E)&\rTo &\bigoplus H^p(E_i,\cO_{E_i})\,.
\end{diagram}
Here the vertical arrows are the canonical maps induced by the inclusion of the
sheaves of locally constant functions into holomorphic functions. As $(V,0)$
has rational singularities the group
$H^p(X,\cO_X)$ vanishes. Thus the map
$\alpha$ in the diagram is the zero map and so the image of $\Gr_p^W
H^p(E,\bbbc)$ under
$\beta$ vanishes. Accordingly $M$ and then $F^p \Gr_p^W H^p(E,\bbbc)$
also vanish, proving (1). 

If $F^q:=F^q H^p(E,\bbbc)$ is the Hodge filtration  and
$W_i:=W_i H^p(E,\bbbc)$ is the weight filtration on
$H^p(E,\bbbc)$ then 
$$
F^p\cap W_{p-1}=0.\leqno(2)
$$
In fact, the Hodge
filtration on each quotient $W_i/W_{i-1}$ is induced by  the
filtration
$F^q\cap W_i$. As $W_i/W_{i-1}$ carries a pure Hodge structure of
weight $i$ it follows that $F^p\cap W_i\subseteq F^p\cap W_{i-1}$ for
all $i<p$ and so $F^p\cap W_{p-1}\subseteq F^p\cap
W_{p-2}\subseteq \cdots \subseteq F^p\cap W_0=0$ vanishes.

From (1) and (2) it follows that 
$$
F^pH^p(E,\bbbc)\cong F^p \Gr_p^W H^p(E,\bbbc)=0. 
$$
By \ref{1.2} (1) $F^pH^p(E,\bbbc)\cong
H^0(E,\Omega^p_X/\Omega_X^p\langle E\rangle (-E))$ and so the latter
module vanishes as well, proving our result.
\qed

\begin{remark}\label{1.4}
We note that \ref{1.1} is no longer true for non-isolated
singularities as is already seen from the case of a product
$V\times\bbbc$. 
\end{remark}

Theorem \ref{1.1} has the following interesting consequence which will
be used later.

\begin{proposition}\label{1.5}
If $(V,0)$ is a rational isolated singularity of dimension $n\ge
2$ then any closed holomorphic p-form $\eta$ on
$V\backslash\{0\}$ with $1\le p\le 2$  is exact, i.e.\ after
shrinking $V$ as a neighbourhood of $0$ we can find a
$(p-1)$-form $\xi$ on
$V\backslash\{0\}$ with $d(\xi)=\eta$.
\end{proposition}

In order to prove this result we need a few preparations. For an
arbitrary complex space $X$ let $H_{DR}^p(X)$ denote the $p$-th
``naive" de Rham cohomology group of $X$
$$
H_{DR}^p(X):=H^p(\Gamma(X,\Omega^\sbullet_{X/\bbbc})).
$$
If $X$ is smooth then $H^p(X,\bbbc)$ is the $p$-th hypercohomology
of the de Rham complex $\Omega^\sbullet_{X/\bbbc}$ and so there is
always a natural map 
$$
\alpha_p:H_{DR}^p(X)\to H^p(X,\Omega^\sbullet_{X/\bbbc})\cong
H^p(X,\bbbc).
$$
Using the the spectral sequence
$E_1^{pq}=H^q(X,\Omega^p_X)\Rightarrow H^{p+q}(X,\bbbc)$
it follows that the map $\alpha_0$ is always bijective and that
$\alpha_1$ is injective; note that the spectral sequence induces an 
exact sequence of small order terms
$$
0\to H^1_{DR}(X)\to H^1(X,\bbbc)\to H^1(X,\cO_X) \to H^2_{DR}(X).
$$
We need the following simple observation.

\begin{lemma}\label{1.6}
Assume that $X$ is smooth. If $H^1(X,\cO_X)=0$ then the natural map
$H_{DR}^1(X)\to H^1(X,\bbbc)$ is bijective and $H_{DR}^2(X)\to
H^2(X,\bbbc)$ is injective.
\end{lemma}

\proof
The first part follows immediately from the sequence of small
order terms above. The second part is easily seen again from the
spectral sequence; we leave the simple details to the reader. 
\qed

\noindent{\it Proof of \ref{1.5}.}
We may assume that $V$ is a contractible
Stein space with 0 as its only singular point. As before let
$\pi:X\to V$ be a resolution of singularities  with exceptional set
$E:=\pi^{-1}(0)$. 
Suppose that $\eta$ is a $p$-form on $U:=V\backslash\{0\}$, where
$p=1$ or $p=2$. By \ref{1.1} the form $\pi^*(\eta)$ extends to a
holomorphic differential form on $X$ that lies in
$\Gamma(X,\Omega_X^p\langle
E\rangle (-E))$. As $\eta$ is closed we can consider its
associated cohomology class $[\pi^*(\eta)]\in H_{DR}^p(X)$. Its image
in
$H^p(X,\bbbc)\cong H^p(E,\bbbc)$ is already contained in 
$F^p(H^p(E,\bbbc))\subseteq H^p(E,\bbbc)$. Under the isomorphism
$$
F^p(H^p(E,\bbbc))\cong H^0(E,\Omega_X^p\langle E\rangle/
\Omega_X^p\langle E\rangle (-E)) 
$$
(see \ref{1.2} (1)) it is represented by the residue
class of $\pi^*(\eta)$  in  $\Omega_X^p\langle E\rangle/
\Omega_X^p\langle E\rangle (-E)$ and so it is zero. Using
\ref{1.6} it follows that the cohomology class of $[\pi^*(\eta)]$ in
$H^p_{DR}(X,\bbbc)$ vanishes. Thus we can find a form $\xi\in
\Gamma(X,\Omega_X^{p-1})$ with $d(\xi)=\eta$, as required.
\qed

\section{Contact singularities}

Recall that a contact form on a complex manifold $X$ of dimension
$2r+1$ is a holomorphic 1-form $\eta$ on $X$ such that 
$$
\eta\wedge (d\eta)^r
$$
has no zero. By Darboux's lemma, in suitable coordinates
$(z, x_1, \ldots, x_{2r})$ such a form can be written as
$$
\eta=dz+\sum_{\varrho=1}^{r} x_{2\varrho-1} dx_{2\varrho}.
$$
By a {\em contact singularity} we mean
a normal singularity
$(V,0)$ together with a 1-form $\eta\in
(\Omega^1_V)^{\vee\vee}_0$ which is a contact form on the regular
part of $V$. Note that the dimension of $V$is necessarily odd so that
we can write
$$
\dim V=2r+1\ge 3.
$$
In analogy with the case of symplectic singularities \cite{Bea} we
have the following facts.

\begin{lemma}\label{c.1}
If $(V,0)$ is a contact singularity then the following hold.

(a) $(V,0)$ is  quasi-Gorenstein, i.e.\ $\omega_{V,0}\cong \cO_{V,0}$.

(b) If $\codim \Sing V\ge 3$ then 
$(V,0)$ is rational. 
\end{lemma}

In particular it follows that a contact singularity with $\codim
\Sing V\ge 3$ is always Gorenstein.

\proof
As $\omega_V$ is generated on $V_\reg$ by the form
$\eta\wedge(d(\eta))^r$, the singularity is quasi-Gorenstein. If
$\codim \Sing V\ge 3$ and $\pi: X\to V$ is a resolution of
singularities, then by \cite{Fle2} the form
$\pi^*(\eta)$ extends to a holomorphic 1-form on $X$ and so
$\pi^*(\eta\wedge(d(\eta))^r)$ extends to a holomorphic form in
$\omega_X$. Using \cite[1.3]{Fle1} it follows that $(V,0)$ is
rational.
\qed

Using \ref{1.5} we can show that the product of a symplectic
singularity with an affine line has the structure of a contact
singularity. Recall that a normal singularity $(V,0)$ of dimension
$2r$ is said to be {\em symplectic} (see \cite[1.1]{Bea}) if
there is a closed 2-form $\eta$ on $U:=\Reg V$ such that
$(d\eta)^{2r}$ has no zero on $U$ near 0 and such that for a
resolution of singularities $\pi:X\to V$ the form
$\pi^*(\eta)$ extends to a holomorphic 2-form on $X$. Note that by
loc.cit.\ such a singularity is always rational.

\begin{proposition}\label{c.2}
Let $(V,0)$ be an isolated
symplectic singularity of dimension
$2r
\ge 2$. Then $(V\times\bbbc, 0)$ is a contact singularity. 
\end{proposition}

\proof
Let $\eta$ be a symplectic form on $V\backslash\{0\}$ so that $\eta$
is a closed non-degenerate holomorphic 2-form. By \ref{1.5} we can
find a holomorphic $1$ form $\xi$ on $V\backslash\{0\}$ with
$d(\xi)=\eta$. With $z$ the local coordinate function on $\bbbc$, the
form 
$$
\omega:= dz+\xi
$$
on $(V\backslash\{0\})\times\bbbc$
then satisfies $\omega\wedge (d\omega)^r=dz\wedge \eta^r$ and so it is
non-degenerate everywhere. Thus $\omega$ is a contact form on
the singularity $(V\times\bbbc,0)$, as required.
\qed

\begin{proposition}\label{c.3}
A complete intersection $(V,0)$ with $\codim \Sing V\ge 3$ can never
carry a contact structure.
\end{proposition}

\proof
In the case of a complete intersection singularity there is an exact
sequence
$$
0\to \cO_V^a\to \cO_V^b\to \Omega^1_V\to 0,
$$
where $b$ is the embedding dimension of $(V,0)$ and $a$ is the number
of defining equations. As $\cO_{V,0}$ is Cohen-Macaulay and
$\codim\Sing X\ge 3$ it follows easily from the Auslander-Buchsbaum
theorem (see e.g.\ \cite{Eis}) that $\Omega^1_{V,0}$ is a reflexive
module. Hence
$\omega$ is a section in $\Omega^1_{V,0}$ and so
$\omega\wedge(d(\omega))^r$ is a section in 
$$
\image(\Omega^n_{V,0}\to \omega_{V,0}).
$$
In particular the cokernel of  $\Omega^n_{V,0}\to
\omega_{V,0}$ is 0 and so we can find $f_1,\ldots, f_{2r+1}$
in the maximal ideal of $\cO_{V,0}$ such that
the form $df_1\wedge\ldots\wedge df_{2r+1}$ generates $\omega_{V,0}$.
Choosing these elements generically we can achieve that the map
$f:=(f_1,\ldots,f_{2r+1}):(V,0)\to(\bbbc^{2r+1},0)$ is finite. By
construction this map is unramified in codimension 1. Using the
theorem on purity of the branch locus (see, e.g.\ \cite[3.2.14]{FOV})
it follows that $f$ is unramified everywhere and so $f$ is an
isomorphism.
\qed

\begin{remark}\label{c.3r}
The same argument shows that $\Omega^1_{V,0}$ cannot be reflexive if
$(V,0)$ is not smooth and carries a contact structure.
\end{remark}

An interesting application --and actually motivation-- of our main
theorem \ref{1.1} is the following result.

\begin{theorem}\label{c.4}
A non-smooth isolated singularity can never carry a contact structure.
\end{theorem}

\proof
Assume that $(V,0)$ is an isolated singularity that carries a
contact structure which is given by the 1-form $\eta$ in
$(\Omega^1_V)^{\vee\vee}$. By definition, $(V,0)$ is then normal of
dimension $\ge 3$. The form
$\delta:= (d\eta)^r$ on $U:=V\backslash\{0\}$ can be considered as a
vector field on $U$ via the canonical isomorphism
$$
\Omega^{2r}_U\cong \Theta_U\otimes \omega_U \cong \Theta_U\;;
$$
note that $\omega_U\cong\cO_U$ by \ref{c.1}. Such a vector field
gives rise to a derivation $\delta:A\to A$, where $A:=\cO_{V,0}$. As
the singularity is isolated $\delta$ maps the maximal ideal of $A$
into itself. Using equivariant resolution of singularities (see e.g.,
\cite{BMi}) we can find a resolution $\pi:X\to V$ such that

(a) $E:=\pi^{-1}(0)$ is a simple normal crossing divisor in $X$ and 

(b) $\delta$ extends to a vector field on $X$ that is tangent to $E$.

\noindent
Using \ref{c.1} the singularity is Gorenstein and rational. By
\ref{1.1} the form $\pi^*(\eta)$ extends to a holomorphic form in
$\Omega^1_X\langle E\rangle (-E)$. As by construction
$\delta$ can be considered as a section in $\Theta_X\langle
E\rangle$ it follows that the contraction
$$
\langle \delta, \pi^*(\eta)\rangle
$$
is a section in $\cO_X(-E)$. Hence on $V$ the contraction $\langle
\delta,\eta\rangle$ must lie in the maximal ideal of $A$ and so
its zero set has codimension 1 in $V$. Under the isomorphism
$\cO_V\cong\omega_V$ the element $\langle
\delta,\eta\rangle$ corresponds to
$\eta\wedge(d\eta)^r$. Thus
$\eta\wedge(d\eta)^r$ has a zero set of codimension 1 as well, which
gives a contradiction.
\qed

We do not know whether there are contact singularities with a
singularity set of even dimension. We even do not know any example of
a contact singularity which is not a product $\bbbc\times V$, where
$V$ is a symplectic singularity.

\section{Extensions of vector fields}

\begin{sit}\label{vf.1}
In this section we consider a non-smooth normal isolated
singularity  $(V,0)$. Let $\pi: X\to V$ be the blowing up of $0\in
V$. We will always assume that $X$ and the exceptional set
$E:=\pi^{-1}(0)$ are smooth.
\end{sit}

For the use in the next section we need the following result.

\begin{proposition}\label{vf.2}
The natural inclusion
$$
H^0(X,(\Lambda^p \Theta_X\langle E\rangle)((p-1)E))\hto
H^0(X,(\Lambda^p \Theta_X\langle E\rangle)(kE))
$$
is an isomorphism for all $k\ge p-1$.
\end{proposition}

For the proof we need the following simple observation.

\begin{lemma}\label{vf.3}
With $(V,0)$ as in \ref{vf.1} and $A:=\cO_{V,0}$ the following hold.

\smallskip

1. Every homomorphism $f:\Omega_A^p\to A$ factors through the maximal
ideal $\fm_A$ of $A$.

2. The image of the natural map $\pi^*(\Omega^p_V)\to \Omega_X^p$ is
just $\Omega_X^p\langle E\rangle(-pE)$. In particular 
$$
\pi^*(\Omega^p_V)/Torsion\cong\Omega_X^p\langle E\rangle(-pE) .
$$
\end{lemma}

\proof
First consider the case $p=1$ in (1). A homomorphism $f:\Omega_A^1\to
A$ corresponds to a derivation $\delta:A\to A$, and it is well known
that such a derivation maps $A$ into $\fm_A$. To deduce (1) for $p>1$
assume that there is a surjective map $f:\Omega_A^p\to A$. For a
sufficiently general form $\eta\in\Omega_A^{p-1}$ the composed map
$\Omega^1_A\stackrel{\wedge \eta}{\lto}\Omega^p_A\to A$ is as well
surjective which contradicts the first part of the proof.

To deduce (2) it is sufficient to treat the case $p=1$; the general
case then follows by taking exterior powers. In the case that
$V=\bbbc^n$ this is easily verified by an explicit computation which
we leave to the reader. In the general case let
$(V,0)\hto(\bbbc^n,0)$ be an embedding. The blowing up $X$ of $V$ in
$0$ then embeds into the blowing up, say, $Y$ of $\bbbc^n$ in $0$,
and $E$ is the intersection of the exceptional set, say $F$, of $Y$
with $X$. Hence the assertion follows from the commutative diagram
with surjective vertical arrows
\begin{diagram}[s=7mm]
\pi^*(\Omega^1_{\bbbc^n})&\rTo & \Omega_Y^1\langle F\rangle(-F)\\
\dTo & & \dTo \\
\pi^*(\Omega^1_V) &\rTo & \Omega_X^1\langle E\rangle(-E).
\end{diagram}
\qed

\proof[Proof of \ref{vf.2}]
A section $\delta$ in $H^0(X,\Lambda^p(\Theta_X\langle E\rangle)(kE))$
gives rise to a map 
$$
\delta: \Omega_X^p\langle E\rangle(-pE)\to
\cO_X((k-p)E).
$$
Composing $\pi_*(\delta)$ with the natural map 
$\Omega^p_V\to\pi_*(\Omega_X^p\langle E\rangle(-pE))$ we get a map
$\Omega^p_V\to\cO_V$, which by \ref{vf.3} (1) factors through the
maximal ideal $\fm_A$ at 0. As by \ref{vf.3} (2)
$\pi^*(\Omega_V^p)\to\Omega_X^p\langle E\rangle(-pE)$ is surjective
the image of $\delta$ is contained in  the ideal sheaf 
$\cO_X(-E)=\fm_A\cO_V$ of $E$ in $\cO_X$. Hence  $\delta$ is in
$H^0(X,(\Lambda^p \Theta_X\langle E\rangle)((p-1)E))$, as required.
\qed

\section{Deformations of symplectic singularities}

In this section let $(V, 0)$ be an isolated symplectic singularity
with $\dim V\ge 4$. We will always assume that the tangent cone $C_0V$
has an isolated singularity. The blowing up $\pi: X\to V$ of
$V$ in $0$ then is smooth with exceptional set $\pi^{-1}(0)$ which is
isomorphic to the projective tangent cone $E:=\bbbp (C_0V)$. Such
singularities were completely classified by \cite{Bea}. We wish to
supplement his result with the following observation.

\begin{proposition}\label{r.1}
 The dimension of the space of infinitesimal
deformations $ T^1_{V,0}$ is given by $b_2(E)-1$.
\end{proposition}

Using the classification in \cite{Bea} it follows that $(V,0)$ is
rigid except when $V$ is the set of $n\times n$-matrices of rank $\le
1$ and of trace 0. 

The proof of this result will follow from
\ref{r.3}, \ref{r.5} below.

Recall that a symplectic singularity of dimension $\ge 4$ is always
rational and Gorenstein \cite{Bea}. We need the following
important facts shown essentially in loc.cit.

\begin{proposition}\label{r.2}
Let $(V,0)$ be an isolated non-smooth symplectic singularity of
dimension
$2r$ with symplectic form $\eta\in ( \Omega^2_{V,0})^{\vee\vee}$. If
the blowing up $\pi:X\to V$ of $0\in V$ is smooth then the following
hold.
\begin{enumerate}
\item $\omega_X\cong \cO_X((r-1)E)$ and
$\omega_E\cong\cO_E(rE)$. 

\item  Multiplying with
$\eta^{r-1}$ gives an isomorphism of vector
bundles
$$
\eta^{r-1}: \Omega_X^1\langle E\rangle (-E)
\lto \Theta_X\langle E\rangle 
\leqno \quad(*)
$$

\item $E$ is a rational homogeneous manifold.

\item $\xi:=\Res \eta\in H^0(E,\Omega^1_E(-E))$ defines a contact
structure on $E$. 

\item $(V,0)$ is homogeneous, i.e.\ there is a locally closed
embedding of
$V\hto \bbbc^N$ such that the ideal of $V$ is given by
homogeneous polynomials.
\end{enumerate}
\end{proposition}

\proof
For the convenience of the reader we include simplified arguments
for some of these facts.  The second part of (1) follows from the
adjunction formula. As the form $\eta^r$ generates $\omega_X$ outside
of $E$, we have
$\omega_X\cong \cO_X(kE)$ for some $k\ge 0$. By \ref{1.1} the form
$\eta$ extends to a section in $\Omega^2_X\langle E\rangle(-E)$
whence $\eta^r$ yields a section in $\Omega^{2r}_X\langle E\rangle
(-rE)\cong \omega_X(-(r-1)E)$. This shows that $k\ge r-1$. 

To show the converse inequality let us first derive (2). As by
\ref{vf.3} (2)
$\Omega^1_X\langle E\rangle (-E)
$ is globally generated  the map $(*)$ factors through
$\Theta\langle E\rangle$ by \ref{vf.2} (applied to the case $p=1$).
Hence we obtain a map as in 
$(*)$ which is an isomorphism outside $E$. Taking the determinant of
this map gives an inclusion
$$
\cO_X((k+1)E-2rE)\cong \det \Omega_X^1\langle E\rangle (-E)
\hto \Theta_X\langle E\rangle \cong \cO_X(-(k+1)E)
$$
and so $k+1-2r\le -(k+1)$, i.e.\ $k+1\le r$. Together with the
inequality $k\ge r-1$ this proves (1). Moreover (2) follows as the
determinant of $(*)$ is an isomorphism.

To deduce (3) note that by (2) and \ref{vf.3}  $\Theta_E$
being a quotient of $\Theta_X\langle E\rangle$ is globally generated
and so $E$ is a rational homogeneous manifold.
(4) is a consequence of the equation $\Res\eta^r=d(\xi)\wedge
\xi^{r-1}$ in $\Omega_E^{2r-1}(-rE)$ which is easily verified by a
local computation using (1). For the proof of (5) we refer the reader
to \cite{Bea}.
\qed

\begin{proposition}\label{r.3}
Let $(V,0)$ be an isolated symplectic singularity of dimension
$n=2r\ge 4$ and let $\pi:X\to V$ be a resolution of singularities such
that
$E:=\pi^{-1}(0)$ is a simple normal crossing divisor. 
Then there is an isomorphism
$$
T^1_{V,0}\cong H^1(X,\Omega^1_X\langle E\rangle). 
$$
\end{proposition}

\proof
We may assume that $V$ is a Stein space and
that  
$U:=V\backslash\{0\}$ is non-singular. We will show that 

(a) $T^1_{V,0}\cong H^1(U,\Theta_U)$, and

(b) $H^1(U, \Theta_U)\cong H^1(X, \Omega^1_X\langle E\rangle)$.

\noindent
To prove (a) we may assume that there is a closed embedding
$V\subseteq W$, where $W$ is a Stein
neighbourhood of 0 in $\bbbc^N$. 
The sheaf
$\Omega^1_W|V\cong \cO_V^N$ is free and so 
$$
\Ext^1_V(\Omega^1_W|V,\cO_V)=0.\leqno (1)
$$
As $(V,0)$ is symplectic, $\cO_{V,0}$ is  rational and Gorenstein (see
\cite[1.3]{Bea}). In particular,
$\depth(\cO_{V,0})=\dim V \ge 3$ and so by \cite{Sche}
$$
H^1(U, \Theta_W|V)\cong H^2_{\{0\}}(V,\cO_V^N)=0.\leqno (2)
$$ 
There is a natural exact sequence 
$$
0\to \cI/\cI^2\to \Omega^1_W|V \to \Omega^1_V\to 0.
$$
Taking the associated cohomology sequences with respect to the
functors 
$$
 \Hom_V(-, \cO_V)\quad\mbox{and}\quad
\Hom_U(-,\cO_U)
$$
and using (1) and (2) we get a commutative diagram
\begin{diagram}[h=8mm]
\Hom(\Omega^1_W|V,\cO_V) & \rTo& \Hom(\cI/\cI^2,\cO_V)&\rTo
&\Ext^1(\Omega^1_V,
\cO_V) &\rTo& 0 \\
\dTo>\cong &&\dTo>\cong &&\dTo>\beta\\
H^0(U,\Theta^1_W|V)&\rTo & H^0(U,(\cI/\cI^2)^\vee) & \rTo
&H^1(U,\Theta_U) &\rTo & 0 \; ,
\end{diagram}
where we have used the fact that for 
locally free modules $\cM$ on $U$ there are isomorphisms 
$H^i(U,\cM^\vee)\cong \Ext^i_U(\cM,\cO_U)$ for all $i$. 
Thus the map $\beta$ in the diagram above is an isomorphism. As
$T^1_V$ is isomorphic to $\Ext^1_V(\Omega^1_V,\cO_V)$, (a) follows.

To deduce (b) note first that the symplectic structure provides an
isomorphism $\Theta_U\cong \Omega^1_U$ and so 
$$
H^1(U, \Theta_U)\cong H^1(U, \Omega^1_U)\leqno (3)
$$
Moreover, by the vanishing theorem for isolated singularities of
Steenbrink \cite{Ste2}
$$
H^q(X,\Omega^p_X\langle E\rangle(-E))=0\quad \mbox{for } p+q>n.
\leqno (4)
$$
By duality this is equivalent to
$$
H^q_E(X,\Omega^p_X\langle E\rangle)=0\quad \mbox{for } p+q<n,
\leqno (4)'
$$
where $H^q_E(...)$ denotes cohomology with support in $E$.
As $U\cong X\backslash E$ there is an exact sequence
$$
H^1_E(X,\Omega^1_X\langle E\rangle)\to
H^1(X,\Omega^1_X\langle E\rangle)\to H^1(U,\Omega^1_U)
\to H^2_E(X,\Omega^1_X\langle E\rangle)\, ,
$$
in which the outer terms vanish because of (2)'. Together with (1)
this proves (b).
\qed

\begin{remark}\label{r.4}
The proof shows more generally that for any isolated singularity
of dimension $\ge 4$ there is always an isomorphism
$H^1(U,\Omega^1_U)\cong H^1(X,\Omega^1_X\langle E\rangle).
$
\end{remark}

Now the remaining argument for the proof of \ref{r.1} is provided by
the following lemma.

\begin{lemma}\label{r.5}
Let $(V,0)$ be as in \ref{r.1} and let $\pi:X\to V$ be the blowing
up of $V$ in 0. Then\/
$\dim_\bbbc H^1(X,\Omega^1_X\langle E\rangle )=b_2(E)-1$.
\end{lemma}

\proof
As the singularity is homogeneous, $\cO_E(1)\cong\cO_X(-E)|E$ and $X$
is just the projective analytic space, say  $\Projan
(\bbbs(\cO_E(1)))$, associated to the symmetric algebra
$\bbbs(\cO_E(1))=\bigoplus_{q\ge 0}\cO_E(q)$. Moreover an easy local
computation shows that 
$\Omega^1_X\langle E\rangle$ is the coherent $\cO_X$-module associated
to the graded $\bbbs(\cO_E(1))$-module $\bigoplus_{q\ge
0}\Omega^1_X\langle E\rangle (-qE)|E)$. This implies that 
$$
H^p(X,\Omega^1_X\langle E\rangle)
\cong
\bigoplus_{q\ge 0}
H^p(E,\Omega^1_X\langle E\rangle (-qE)|E).
$$
It suffices to show that 

(a) $H^p(E,\Omega^1_X\langle E\rangle (-qE)|E)=0$ for $p> 0$ and
$q>0$.

(b) $\dim_\bbbc H^1(E,\Omega^1_X\langle E\rangle |E)=b_2(E)-1$.

\noindent
To deduce (a) note first that by \ref{r.6} below the groups $H^p(E,
\Theta_E(-qE))$ vanish if $p>0$ and
$q\ge 0$. Using the exact sequence 
$$
0\lto \cO_E(-qE)\lto \Theta_X\langle E\rangle (-qE)|E \lto 
\Theta_E(-qE) \lto 0
$$
it follows that $H^p(E,\Theta_X\langle E\rangle (-qE)|E)$
also vanishes if $p>0$ and $q\ge 0$. By \ref{r.2} (2) this gives
$$
H^p(E,\Omega^1_X\langle E\rangle (-(q+1)E)|E)\cong
H^p(E,\Theta_X\langle E\rangle (-qE)|E) =0
$$
for $p>0$ and $q\ge 0$,
proving (a). To deduce (b), let us consider
the exact sequence 
$0\to \Omega^1_E\to \Omega^1_X\langle E\rangle |E \to 
\cO_E \to 0$ and its cohomology sequence
$$
\bbbc\cong H^0(E,\cO_E)\stackrel{\partial}{\lto} H^1(E,
\Omega^1_E)\to  H^1(E,
\Omega^1_X\langle E\rangle|E) \to H^1(E, \cO_E)=0.
$$
The map $\partial$ is multiplication with the Chern class of
$\cO_E(E)$ and so is injective. This immediately gives (b).
\qed

In the proof above we have used the following well known fact for
which we could not find a reference.

\begin{lemma}\label{r.6}
Let $X=G/P$ be a rational homogeneous manifold, where $G$ is a
semisimple linear algebraic group and $P$ is a parabolic subgroup of
$G$. If $\cL$ is an ample line bundle on $X$ then
$H^p(X,\Theta_X\otimes_{\cO_X}\cL^{\otimes q})=0$ for all $p\ge 1$
and $q\ge 0$. 
\end{lemma}

\proof
This follows with the same arguments as in \cite[4.8 Theorem 1]{Akh}
as the highest weight of an ample line bundle on $X$ is known to be
positive. 
\qed

\begin{remark}
After completing this paper we were informed that Lemma \ref{1.3} was
independently shown in a recent preprint by Y.~Namikawa (cf.\
\cite[1.2]{Na}). 
\end{remark}

\end{document}